\documentclass[11pt]{amsart}

\usepackage[parfill]{parskip}   
\usepackage{amsfonts}
\usepackage{amsmath}
\usepackage{graphicx}
\usepackage{epstopdf}

\newtheorem{theorem}{Theorem}
\newtheorem{lemma}{Lemma}
\newtheorem{corollary}{Corollary}

\title{A unimodular bijection between harmonic vectors of 2-isomorphic graphs} 

\author{William Watkins\\California State University, Northridge}
\thanks{bill.watkins@csun.edu}

\newcommand{\Har}{\text{Harm}}
\newcommand{\Ut}{U_\tau}
\newcommand{\As}{\mathcal{A}}
\newcommand{\Cs}{\mathcal{C}}
\newcommand{\Ss}{\mathcal{S}}
\newcommand{\Ng}{N_G}
\newcommand{\Nh}{N_H}
\newcommand{\RS}{\text{rowspace}}

\begin{document}

\date{June 12, 2026}

\maketitle

\begin{abstract}
Let $G$ and $H$ be connected graphs that are 2-isomorphic.  It is known that their Laplacian matrices are congruent by a unimodular matrix $U$. 
 In this paper we show (Thm. \ref{thm:main2}) 
 that $U$ is a bijection between  certain spaces of harmonic vectors on the vertices  of $G$ and $H$.
In particular  (Cor. \ref{cor:main1}) 
if $u$ is a harmonic vector with respect to vertices $c, d$ in $H$ and the 2-isomorphism maps edge $(a,b)$ in $G$ to edge $(c,d)$ in $H$, then $uU$ is a harmonic vector with respect to vertices $a, b$ in $G$.  
\end{abstract}

{\bf Key words:}
{graph,  2-isomorphism, Laplacian matrix,  unimodular congruence, harmonic vector}

{\bf AMS subject classification:}
 05C22, 05C50, 15A63, 15B99

\section{Introduction} Let $G$ and $H$ be connected graphs with $n$ edges and Laplacian matrices $L_G, L_H$.  
It is known that $G$ and $H$ are 2-isomorphic if and only 
there exists a unimodular matrix $U$ such that $UL_GU^T=L_H$. 
The $n \times n$ unimodular matrix sends row vectors indexed by the vertices of $G$ onto row vectors indexed by the vertices of $H$. Thus $U$ provides a strong link between these row vectors.
And since (vertex) harmonic vectors are defined in terms of the Laplacian matrix, it is reasonable to suspect that $U$ may send harmonic vectors in $G$ to harmonic vectors in $H$.    In this paper, we show that this is true.

Harmonic vectors play an important role in the theory of absorbing Markov chains and in the theory of electrical resistor networks.  See \cite{DS} and the groundbreaking work of Foster and Whitney \cite{Fos,Whi1,Whi2}.

Definitions and preliminary results are given Section \ref{sec:prelim}.  The main results are stated in Section \ref{sec:main} along with an example. Proofs  are in Section \ref{sec:proofs} and a generalization to graphs with weighted edges is in Section \ref{sec:addendum}.

\section{Preliminary definitions and results} \label{sec:prelim}
\subsection{Incidence and Laplacian matrices of a graph}
Let $G$ be a connected graph with $n$ vertices $\{1, 2, \ldots, n\}$ and $m$ edges $E(G)=\{ (a_1, b_1), \ldots, (a_m, b_m)\}$.  To each edge $(a_l,b_l)$ assign a positive end and a negative end, say $a_l$ is the positive end and $b_l$ is the negative end.  The {\it vertex-edge incidence matrix}, $\Ng$, corresponding to this order and plus/minus-assignments is the $n \times m$ matrix defined as follows:
Column $l$ has only two nonzero entries; $+1$ in position $a_l$ and  $-1$ in position $b_l$.  
The {\it Laplacian matrix} is defined by $L_G=\Ng\Ng^T$.  The Laplacian matrix does not depend on the order or the plus/minus assignments in the columns of $\Ng$.  Indeed, $L_G$ can be defined without using an incidence matrix:
\[
{(L_G)}_{a,b} = 
\begin{cases}
-1& \text{ if } a\neq b \text{ and } (a,b) \in E(G)\\
\ \ 0 & \text{ if } a\neq b \text{ and } (a,b) \notin E(G)\\
\ \ d_a &\text { if } a=b,
\end{cases}
\]
where $d_a$ is the degree of vertex $a$.
See \cite{Mer} for more on the Laplacian matrix.

\subsection{2-isomorphism} \label{sec:cycleiso}
Let $H$ be another connected graph with $n$ edges $\{1, 2, \ldots, n\}$ and $m$ edges $E(H)$.  A bijection $\tau: E(G) \rightarrow E(H)$ is a {\it 2-isomorphism} \cite{Wel, Ox} if a subset $S \subseteq E(G)$ of edges form a cycle in $G$ if and only if $\tau(S)$ form a cycle in $H$.  In that case $G$ and $H$ are said to be {\it 2-isomorphic}.  (See the seminal work on 2-isomorphism by Whitney \cite{Whi1, Whi2}.)

Of interest here is a condition on the Laplacian matrices, $L_G, L_H$, and on the incident matrices, $N_G, N_H$, that are equivalent to $G$ and $H$ being 2-isomorphic.

\begin{theorem}[\cite{Wat1,Wat2}] \label{thm:Wat}
Let $G$ and $H$ be graphs with $n$ vertices.  Then $G$ and $H$ are 2-isomorphic if and only if
there exist a unimodular matrix $U$ such that $L_H=UL_GU^T$.  Furthermore if $L_H=UL_GU^T$ for a unimodular matrix $U$, then $U\Ng$ is an incidence matrix for $H$ whenever $\Ng$ is an incidence matrix for $G$.
\end{theorem}

\subsection{Harmonic vectors} 
 Let $\As \subseteq \{1, \ldots, n\}$ be a set of $k$ vertices in $G$.  A row vector $v=(v_1 \ldots, v_n)$ is {\it harmonic} with respect to $\As$  if 
 every vertex $j\notin \As$ satisfies an averaging condition:
 $v_j$ is the average of its neighbors.  That is
 \[
 v_j=\frac{1}{d_j}\sum_{k\in \text{Nbhd}(j)} v_k,
 \]
 where $d_j$ is the degree of vertex $j$ and the sum is taken over all neighbors of vertex $j$.
 This is precisely the condition that
 the $j$th coordinate of $vL_G$ is zero for all $j \notin \As$.  In addition, the coordinates of $v L_G$ must sum to zero since $L_G e^T=0$.  
Thus we  define (and denote) the harmonic vectors of $G$, with respect to vertices $\As$, by
\[
\Har(G,\As) = \{v \in \mathbb{R}^n: vL_G \in W(\As)\},
\]
where
\[
W(\As)=\{u=(u_1, \ldots, u_n): u_i=0 \text{ for } i\notin \As 
\text{ and } \sum_{i=1}^n u_i=0\}.
\]

We record some simple facts about harmonic vectors in the next Lemma:
\begin{lemma} \label{lem:harmonic}
Let $G$ be a connected graph,  and let $\As$ be a subset of $k$ vertices in $G$.  Then \begin{enumerate}
\item $\dim \Har(G,\As)=k$
\item the vector $e=(1,1, \dots. 1) \in \Har(G,\As)$
\item if $\{\alpha_i: i \in \As\}$ are real numbers, there exists exactly one vector $(v_1, \ldots, v_n) \in \Har(G,\As)$ with $v_i=\alpha_i$ for all $i \in \As$.
\end{enumerate}
\end{lemma}

Part (3) of the lemma tells us that a harmonic vector is determined by its coordinates $v_i$ for $i \in \As$.  Indeed a {\it canonical} basis $b_1, \ldots, b_k$, with $b_i=(b_{i1},\ldots, b_{in})$, for $\Har(G,\As)$ can be chosen so that for $i,j \in \As$
\[
b_{ij}=\begin{cases}
	1, \text{ if }i=j\\
	0, \text{ if }i\neq j.
	\end{cases}
\]

E.g. if $n=6$ and $\As=\{2,4\}$ then the canonical basis is of the form:
\[
b_1= (\ast,1,\ast,0,\ast,\ast), \qquad b_2=(\ast,0,\ast,1,\ast,\ast),
\]
where $\ast$ denotes an unspecified coordinate.

 It follows from the second statement in Theorem \ref{thm:Wat} that if $UL_GU^T=L_H$ and $(a,b) \in E(G)$, then $U(e_a-e_b)^T=(e_c-e_d)^T$ for some $c,d \in \{1, 2, \ldots, n\}$.

However,  it is implicit in \cite{Wat2} that if $\tau: E(G) \rightarrow E(H)$ is a 2-isomorphism, then $U$ can be chosen so that 
$(c,d)=\tau(a,b)$.  If the unimodular matrix $U$ has this property, we say that $U$ is {\it compatible} with $\tau$ and we denote it by $\Ut$.
We summarize this result for reference in the following lemma:
\begin{lemma} \label{lem:compatible}
Let $\tau: E(G) \rightarrow E(H)$ be a 2-isomorphism.  Then there exists a unimodular matrix $\Ut$, compatible with $\tau$, such that $\Ut L_G\Ut^T=L_H$.   That is, if $(a,b) \in E(G)$ then 
\begin{equation} \label{eqn:compatible}
\Ut (e_a-e_b)^T = (e_c-e_d)^T \text{ where } (c,d)=\tau(a,b).
\end{equation}
Equivalently, if $\Ng$ is an incidence matrix for $G$ then the incidence matrix $\Nh=\Ut \Ng$ for $H$  satisfies this property: if Column $l$ of $\Ng$ is $(e_a-e_b)^T$ and $\tau(a,b)=(c,d)$, then Column $l$ of $\Nh$ is $(e_c-e_d)^T$.
\end{lemma}

\section{Main result} \label{sec:main}
The main result shows that if $G$ and $H$ are 2-isomorphic, then there is a unimodular bijection between certain spaces of harmonic vectors on $G$ and on $H$.

\begin{theorem}\label{thm:main2}
Let $G$ and $H$ be connected graphs with $n$ vertices, $\{1,2,\ldots,n\}$, and $m$ edges, and let $\tau: E(G) \rightarrow E(H)$ be a 2-isomorphism.
Let $\Ss_G$ be a subset of $E(G)$ containing $k-1$ edges and let $\Ss_H=\tau(\Ss_G)$ be the corresponding subset of $E(H)$.  Assume that the edges in $S_G$ form a subtree in $G$ and that the edges in $\Ss_H$ form a subtree in $H$.  Let $\As$ be the set of vertices spanned by the edges in $\Ss_G$ and 
$\Cs$ be the set of vertices spanned by the edges in $\Ss_H$.  Then
\[
\Har(G,\As)= \Har(H,\Cs)\Ut.
\]
\end{theorem}
{\bf Note}: A set of edges $\Ss_G$ is a {\it subtree} if it is connected and has no cycles.

Corollary \ref{cor:main1} is the special case of Theorem \ref{thm:main2} where $\Ss_G, \Ss_H$ are subtrees containing a single edge.
\begin{corollary}\label{cor:main1}
Let $G$ and $H$ be connected graphs with $n$ vertices, $\{1,2,\ldots,n\}$, $m$ edges $E(G), E(H)$, and Laplacian matrices $L_G, L_H$. Let $\tau: E(G) \rightarrow E(H)$ be a 2-isomorphism and let 
$\Ut$ be a unimodular matrix such that $\Ut L_G\Ut^T=L_H$ and $\Ut$ is compatible with $\tau$, as guaranteed by Lemma \ref{lem:compatible}.  Suppose  $(a,b)$ is an edge in $G$ and $\tau(a,b)=(c,d)$ is the corresponding edge in $H$.  Then
\[
\Har(G,\{a,b\}) = \Har(H,\{c,d\})\Ut.
\]
\end{corollary}
\subsection{An example}
Before starting the proof, we give a simple example of Corollary \ref{cor:main1}.
The graphs below are 2-isomorphic
\[
\begin{array}{cc}
\includegraphics[width=1.2in]{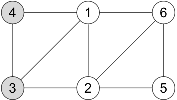}  &\includegraphics[width=1.2in]{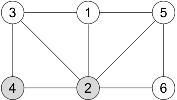} \\
G, \{a,b\}=\{3,4\}&H, \{c,d\}=\{2,4\}
\end{array},
\]
with 2-isomorphism $\tau$:
\[
\begin{array}{r|ccccccccc}
 (a,b)\in E(G)&(1,2) & (1,3) & (1,4) & (1,6) & (2,3) & (2,5) & (2,6) & \mathbf{(3,4)} & (5,6) \\
 \tau(a,b)\in E(H)&(1,2) & (3,2) & (3,4) & (1,5) & (3,1) & (2,6) & (2,5) & \mathbf{(2,4)} & (6,5) \\
\end{array},
\]
and corresponding unimodular matrix
\[
U_{\tau}=\left[
\begin{array}{rrrrrr}
 1 & 0 & 1 & 1 & 0 & 0 \\
 0 & 1 & 1 & 0 & 0 & 0 \\
 0 & 0 & -1 & -1 & 0 & 0 \\
 0 & 0 & 0 & 1 & 0 & 0 \\
 0 & 0 & 0 & 0 & 0 & 1 \\
 0 & 0 & 0 & 0 & 1 & 0 \\
\end{array}
\right].
\]
Take $(a,b)=(3,4)$, $(c,d)=\tau(a,b)=(2,4)$.  
It is easy to check that $\Har(H,\{2,4\})=\RS(B)$ with the matrix $B$ below:
\[
B= \frac{1}{34}\left[
\begin{array}{cccccc}
 5 & 0 & 13 & 34 & 2 & 1 \\
 29 & 34 & 21 & 0 & 32 & 33
\end{array}
\right],
\]
whose rows are the canonical basis for $\Har(H,\{2,4\})$.

Corollary \ref{cor:main1} guarantees that $\Har(G,\{3,4\})=\Har(H,\{2,4\})\Ut$.  But $\Har(H,\{2,4\})\Ut=\RS(B\Ut)$ and
\[ 
B \Ut =\frac{1}{34}\left[
\begin{array}{cccccc}
 5 & 0 & -8 & 26 & 1 & 2 \\
 29 & 34 & 42 & 8 & 33 & 32
\end{array}
\right].
\]
The rows of $B \Ut$ do not form the canonical basis for 
$\Har(G,\{3,4\})$. But it is easy to check that
\[
\RS(B \Ut) =\RS(A)
\]
for the matrix $A$: 
\[
A=\frac{1}{34}\left[
\begin{array}{cccccc}
 21 & 26 & 34 & 0 & 25 & 24 \\
 13 & 8 & 0 & 34 & 9 & 10 
\end{array}
\right],
\]
whose rows do form the canonical basis for $\Har(G,\{3,4\})$.

\section{Proof of Theorem \ref{thm:main2}} \label{sec:proofs}
Let $G$ and $H$ be connected, 2-isomorphic graphs with edge map $\tau: E(G) \rightarrow E(H)$ and corresponding unimodular matrix $\Ut$ such that $\Ut N_G=N_H$.  Also assume that the $k-1$ columns corresponding to $\Ss_G$ and to $\Ss_H$ are listed first in $N_G$ and $N_H$.  Then
\[
N_G=\begin{bmatrix}
	N_1 : N_2
	\end{bmatrix}
\text{ and }
N_H=\begin{bmatrix}
	M_1 :  M_2
	\end{bmatrix},
\]
where $N_1, M_1$ are $(k-1) \times n$ matrices whose columns correspond to the edges in $\Ss_G$, $\Ss_H$, respectively.  
\subsection{$\RS (N_1^T)=W(\As) \text{ and } \RS (M_1^T)=W(\Cs)$}\label{subsec:proof1}
The only nonzero entries in the columns of $N_1$ (the rows of $N_1^T$) occur in positions $\As$ so  $\RS (N_1^T)\subseteq W(\As)$.  Furthermore the $k-1$ columns of $N_1$ are linearly independent because $\Ss_G$ is a subtree.  Therefore the dimension of $\RS(N_1^T)$ is $k-1$.
The dimension of $W(\As)$ is also $k-1$. Therefore 
$\RS (N_1^T)=W(\As)$ and likewise  $\RS (M_1^T)=W(\Cs)$.
\subsection{$W(\As) \Ut^T=W(\Cs)$} \label{subsec:proof2}
From $\Ut N_G=N_H$ we get $\Ut N_1=M_1$ so that $N_1^T \Ut^T=M_1^T$. Therefore 
$\RS(N_1^T \Ut^T)=\RS(M_1^T)$, which combined with \ref{subsec:proof1} gives \ref{subsec:proof2}.

\subsection{$u \in \Har(H,\Cs)$ and $v=u\Ut$ implies $v \in \Har(G,\As)$}
  Assume $u \in \Har(H,\Cs)$ and $v=u\Ut$. Then
\begin{equation} \label{eqn:WtoW}
vL_G \Ut^T=u\Ut L_G \Ut^T=u L_H=y\in W(\Cs).
\end{equation}
From \ref{subsec:proof2} there exists $x \in W(\As)$. such that  $y=x \Ut^T$.
Thus $vL_G \Ut^T=x \Ut$ and $vL_G=x \in W(\As)$. 
So $v \in \Har(G,\As)$.\\
$\square$

\section{Addendum: weighted edges} \label{sec:addendum}
  Theorem \ref{thm:main2} also holds for graphs with weighted edges.  Replace $L_G$ with a   {\it weighted-Laplacian} matrix
$L_G^w$ where 
\[
L_G^w=N_G D N_G^T,
\] 
$D$ is an $m \times m$ diagonal matrix of edge-weights, and
\[
\Har(G,w,\As)=\{v\in \mathbb{R}^n : vL_G^w\in W(\As)\}.
\]
The proofs for the weighted graphs are the same as for the unweighted graphs.


\begin{thebibliography}{99999}

\bibitem {Mer} Merris R. Laplacian matrices of graphs: a survey. Linear 
Algebra  Appl. 1994;197,198:143-176.

\bibitem {Wel} Welsh DJA., Matroid Theory. Academic Press;1976.


\bibitem {Ox} Oxley J. Matroid Theory. 2nd ed. Oxford Graduate Texts in Mathematics. Oxford University Press;2011.

\bibitem {Whi1} Whitney H. On the classification of graphs. American Journal of Mathematics. 1933;55(1):236-244.

\bibitem {Whi2} Whitney H. 2-isomorphic graphs, American Journal of Mathematics. 1933;55(1):245-254.

\bibitem {Wat1}  Watkins W. The Laplacian matrix of a graph: unimodular congruence. Linear Multilinear Algebra. 1990;28:35-43.

\bibitem {Wat2} Watkins W. Unimodular congruence of the Laplacian matrix of a graph. Linear Algebra Appl. 1994;201:43-49.

\bibitem {DS} Doyle PG, Snell JL. Random walks and electrical networks. The Carus Mathematical Monographs No. 22. Mathematical Association of America;1984.

\bibitem {Fos} Foster RM. Geometrical circuits of electrical networks. Transactions A.I.E.E.  June 1932:309-317.



\end{thebibliography}
\end{document}